\documentclass[12pt,a4paper]{amsart}
\usepackage{eepic}

\newcommand{\be}{\begin{equation}}
\newcommand{\ee}{\end{equation}}

\newcommand{\g}{{\mathfrak g}}
\newcommand{\uqgh}{U_q(\hat\g)}

\newcommand{\Rc}{\check{R}}
\newcommand{\C}{\mathbb{C}}
\newcommand{\Q}{\mathbb{Q}}

\newcommand{\RR}{{\mathcal{R}}}
\newcommand{\h}{\mathfrak{h}}

\DeclareMathOperator{\tr}{tr}
\DeclareMathOperator{\End}{End}

\DeclareMathOperator{\id}{id}

\numberwithin{equation}{section}

\newcommand{\drawpath}[4]{\path(#1,#2)(#3,#4)}
\newcommand{\drawthickdot}[2]{\put(#1,#2){\circle*{1}}}%
\newcommand{\drawcenteredtext}[3]{\put(#1,#2){\makebox(0,0){#3}}}
\newcommand{\drawrighttext}[3]{\put(#1,#2){\makebox(0,0)[r]{#3}}}
\newcommand{\drawlefttext}[3]{\put(#1,#2){\makebox(0,0)[l]{#3}}}

\begin{document}

\title{Affine quantum groups}
\author{G.W. Delius and N.J. MacKay}
\address{Department of Mathematics\\
University of York, UK}
\thanks{Contribution to the Encyclopedia of Mathematical Physics, Elsevier 2006}
\maketitle

Affine quantum groups are certain pseudo-quasitriangular Hopf
algebras that arise in mathematical physics in the context of
integrable quantum field theory, integrable quantum spin chains, and
solvable lattice models. They provide the algebraic framework behind
the spectral parameter dependent Yang-Baxter equation
\begin{equation}\label{YB}
R_{12}(u)R_{13}(u+v)R_{23}(v)=R_{23}(v)R_{13}(u+v)R_{12}(u).
\end{equation}

One can distinguish three classes of affine quantum groups, each
leading to a different dependence of the R-matrices on the spectral
parameter $u$: Yangians lead to rational R-matrices, quantum affine
algebras lead to trigonometric R-matrices and elliptic quantum
groups lead to elliptic R-matrices. We will mostly concentrate on
the quantum affine algebras but many results hold similarly for the
other classes.

After giving mathematical details about quantum affine algebras and
Yangians in the first two section, we describe how these algebras
arise in different areas of mathematical physics in the three
following sections. We end with a description of boundary quantum
groups which extend the formalism to the boundary Yang-Baxter
(reflection) equation.

\section{Quantum affine algebras}

\subsection{Definition}

A quantum affine algebra $\uqgh$ is a quantization of the enveloping
algebra $U(\hat\g)$ of an affine Lie algebra (Kac-Moody algebra)
$\hat\g$. So we start by introducing affine Lie algebras and their
enveloping algebras before proceeding to give their quantizations.

Let $\g$ be a semisimple finite-dimensional Lie algebra over $\C$ of
rank $r$ with Cartan matrix $(a_{ij})_{i,j=1,...,r}$, symmetrizable
via positive integers $d_i$, so that $d_i a_{ij}$ is symmetric. In
terms of the simple roots $\alpha_i$, we have
$$
 a_{ij} = 2\frac{\alpha_i\cdot\alpha_j}{|\alpha_i|^2}
 \qquad {\rm and} \qquad d_i
 = {|\alpha_i|^2\over 2}.
$$
We can introduce an $\alpha_0 = \sum_{i=1}^r n_i\alpha_i$ in such
a way that the extended Cartan matrix $(a_{ij})_{i,j=0,\dots,r}$
is of affine type -- that is, it is positive semi-definite of rank
$r$. The integers $n_i$ are referred to as Kac indices. Choosing
$\alpha_0$ to be the highest root of $\g$ leads to an untwisted
affine Kac-Moody algebra while choosing $\alpha_0$ to be the
highest short root of $\g$ leads to a twisted affine Kac-Moody
algebra.

One defines the affine Lie algebra $\hat\g$ corresponding to this
affine Cartan matrix as the Lie algebra (over $\C$) with generators
$H_i, E^\pm_i$ for $i=0,1,...,r$ and $D$ and relations
\begin{gather}\label{lierels}
 [H_i, E^\pm_j]   =  \pm a_{ij} E^\pm_i\,,
 \ \ [H_i,H_j]  =  0\,,\nonumber\\
 [E^+_i , E^-_j]  = \delta_{ij} H_i,\\
 [D,H_i]=0,\ \ \ [D,E^\pm_i]=\pm\delta_{i,0} E^\pm_i,\nonumber\\
 \sum_{k=0}^{1-a_{ij}} (-1)^k \left( \begin{array}{c} 1-a_{ij} \\ k
 \end{array} \right) (E_i^\pm)^k E_j^\pm
 (E_i^\pm)^{1-a_{ij}-k}=0,
 \ \  i\neq j.\nonumber
\end{gather}
The $E^\pm_i$ are referred to as Chevalley generators and the last
set of relations are known as Serre relations. The generator $D$ is
known as the canonical derivation. We will denote the algebra
obtained by dropping the generator $D$ by $\hat\g'$.

In applications to physics the affine Lie algebra $\hat\g$ often
occurs in an isomorphic form as the loop Lie algebra
$\g[z,z^{-1}]\oplus\mathbb{C}\cdot c$ with Lie product (for
untwisted $\hat\g$)
\begin{equation}\label{loop}
 [Xz^k,Yz^l]=[X,Y]z^{k+l}+\delta_{k,-l}(X,Y)c,\ \ \text{for}\ \
 X,Y\in\g,\ k,l\in\mathbb{Z},
\end{equation}
and $c$ being the central element.

The universal enveloping algebra $U(\hat\g)$ of $\hat\g$ is the
unital algebra over $\C$ with generators $H_i, E^\pm_i$ for
$i=0,1,...,r$ and $D$ and with relations given by \eqref{lierels}
where now $[\ ,\ ]$ stands for the commutator instead of the Lie
product.

To define the quantization of $U(\hat\g)$ one can either define
$U_h(\hat\g)$ \cite{D1} as an algebra over the ring $\C[[h]]$ of
formal power series over an indeterminate $h$ or one can define
$U_q(\hat\g)$ \cite{J0} as an algebra over the field $\Q(q)$ of
rational functions of $q$ with coefficients in $\Q$. We will present
$U_h(\hat\g)$ first.

The {\em quantum affine algebra} $U_h(\hat\g)$ is the unital algebra
over $\C[[h]]$ topologically generated by $H_i, E^\pm_i$ for
$i=0,1,...,r$ and $D$ with relations
\begin{gather}\label{qrels}
 [H_i, E^\pm_j]   =  \pm a_{ij} E^\pm_i\,, \ \ [H_i,H_j]  =  0\,,\nonumber\\
 [E^+_i , E^-_j]  = \delta_{ij} \frac{q_i^{H_i}-q_i^{-H_i}}
 {q_i-q_{-i}}\,,\\
 [D,H_i]=0\ \ \ [D,E^\pm_i]=\pm\delta_{i,0} E^\pm_i,\nonumber\\
 \sum_{k=0}^{1-a_{ij}} (-1)^k \left[ \begin{array}{c} 1-a_{ij} \\ k
 \end{array} \right]_{q_i} (E_i^\pm)^k E_j^\pm
 (E_i^\pm)^{1-a_{ij}-k}=0,
 \ \  i\neq j,\nonumber
\end{gather}
where $q_i=q^{d_i}$ and $q=e^h$. The $q$-binomial coefficients are
defined by
\begin{gather}\label{qnums}
 [n]_q=\frac{q^n-q^{-n}}{q-q^{-1}},\\
 [n]_q!=[n]_q\cdot[n-1]_q,\cdots[2]_q[1]_q,\\
 \left[\begin{array}{c} m \\ n \end{array} \right]_{q}
 = \frac{[m]_q!}{[n]_q![m-n]_q!}.
\end{gather}

The quantum affine algebra $U_h(\hat\g)$ is a Hopf algebra with
co-product
\begin{eqnarray}\label{coproduct}
\Delta(D) & = & D \otimes 1 + 1 \otimes D ,\nonumber\\
\Delta(H_i)& = & H_i\otimes 1 + 1\otimes H_i, \\
\Delta(E^\pm_i) & = & E^\pm_i \otimes q_i^{- H_i/2} + q_i^{ H_i/2}
\otimes E^\pm_i,\nonumber
\end{eqnarray}
antipode
\begin{equation}\label{antipode}
 S(D)=-D,\quad S(H_i)=-H_i,\quad
 S(E^\pm_i)=-q_i^{\mp 1} E^\pm_i,
\end{equation}
and co-unit
\begin{equation}\label{counit}
 \epsilon(D)=\epsilon(H_i)=\epsilon(E^\pm_i)=0.
\end{equation}

It is easy to see that the classical enveloping algebra $U(\hat\g)$
can be obtained from the above by setting $h=0$, or more formally,
$$
 U_h(\hat\g)/h U_h(\hat\g) = U(\hat\g).
$$

We can also define the quantum affine algebra $U_q(\hat\g)$ as the
algebra over $\Q(q)$ with generators $K_i, E^\pm_i, D$ for
$i=0,1,...,r$ and relations that are obtained from the ones given
above for $U_h(\hat\g)$ by setting
\begin{equation}\label{Kdef}
 q_i^{H_i/2} = K_i,\ \  i = 0,\dots,r.
\end{equation}
One can also go further to an algebraic formulation over $\C$ in
which $q$ is a complex number (with some points including $q=0$ not
allowed). This has the advantage that it becomes possible to
specialise for example to $q$ a root of unity, where special
phenomena occur.

\subsection{Representations}

For applications in physics the finite-dimensional representations
of $U_h(\hat\g')$ are the most interesting. As will be explained in
later sections, these occur for example as particle multiplets in
2-d quantum field theory or as spin Hilbert spaces in quantum spin
chains. In the next subsection we will use them to derive matrix
solutions to the Yang-Baxter equation.

While for non-affine quantum algebras $U_h(\g)$ the ring of
representations is isomorphic to that of the classical enveloping
algebra $U(\g)$ (because in fact the algebras are isomorphic, as
Drinfeld has pointed out), the corresponding fact is no longer true
for affine quantum groups, except in the case $\hat\g$ = $a_n^{(1)}$
= $\widehat{{\mathfrak{sl}}_{n+1}}$.

For the classical enveloping algebras $U(\hat\g')$ any
finite-dimensional representation of $U(\g)$ also carries a
finite-dimensional representation of $U(\hat\g')$. In the quantum
case however in general an irreducible representation of
$U_h(\hat\g')$ reduces to a sum of representations of $U_h(\g)$.

To classify the finite-dimensional representations of $U_h(\hat\g')$
it is necessary to use a different realization of $U_h(\hat\g')$
that looks more like a quantization of the loop algebra realization
\eqref{loop} than the realization in terms of Chevalley generators.
In terms of the generators in this alternative realization, which we
do not give here because of its complexity, the finite-dimensional
representations can be viewed as pseudo-highest weight
representations. There is a set of $r$ `fundamental' representations
$V^a$, $a=1,...r$, each containing the corresponding $U_h(\g)$
fundamental representation as a component, from the tensor products
of which all the other finite-dimensional representations may be
constructed. The details can be found in \cite{CP}.

Given some representation $\rho: U_h(\hat\g') \rightarrow \End(V)$
we can introduce a parameter $\lambda$ with the help of the
automorphism $\tau_\lambda$ of $U_h(\hat\g')$ generated by $D$ and
given by
\begin{equation}\label{tau}
 \tau_\lambda(E^\pm_i) = \lambda^{\pm s_i} E^\pm_i,\ \
 \tau_\lambda(H_i)=H_i,\ \ i=0,\dots,r.
\end{equation}
Different choices for the $s_i$ correspond to different gradations.
Commonly-used are the {\em homogeneous gradation}, $s_0=1,
s_1=...=s_r=0$, and the {\em principal gradation},
$s_0=s_1=...=s_r=1$. We shall also need the {\em spin gradation}
$s_i=d_i^{-1}$. The representations
$$
 \rho_\lambda = \rho\circ\tau_\lambda
$$
play an important role in applications to integrable models where
$\lambda$ is referred to as the (multiplicative) spectral parameter.
In applications to particle scattering introduced in a later section
it is related to the rapidity of the particle. The generator $D$ can
be realized as an infinitesimal scaling operator on $\lambda$ and
thus plays the role of the Lorentz boost generator.

The tensor product representations $\rho^a_\lambda\otimes\rho^b_\mu$
are irreducible generically but become reducible for certain values
of $\lambda/\mu$, a fact which again is important in applications
(fusion procedure, particle bound states).

\subsection{R-matrices}

A Hopf algebra $A$ is said to be {\em almost cocommutative} if there
exists an invertible element $\RR\in A\otimes A$ such that
\begin{equation}\label{Rint}
 \RR \Delta(x) = (\sigma \circ \Delta(x)) \RR,\ \ \text{for all}\ x\in A,
\end{equation}
where $\sigma:x\otimes y \mapsto y \otimes x$ exchanges the two
factors in the coproduct. In a {\em quasitriangular} Hopf algebra
this element $\RR$ satisfies
\begin{equation}\label{DeltaR}
 (\Delta\otimes\id)(\RR) = \RR_{13}\RR_{23},\ \
 (\id\otimes\Delta)(\RR) = \RR_{13}\RR_{12}
\end{equation}
and is known as the {\em universal R-matrix} [See also Article 28].
As a consequence of \eqref{Rint} and \eqref{DeltaR} it automatically
satisfies the Yang-Baxter equation
\begin{equation}
 \RR_{12}\RR_{13}\RR_{23}=\RR_{23}\RR_{13}\RR_{12}.
\end{equation}

For technical reasons, to do with the infinite number of root
vectors of $\hat\g$, the quantum affine algebra $U_h(\hat\g)$ does
not possess a universal R-matrix that is an element of
$U_h(\hat\g)\otimes U_h(\hat\g)$. However, as pointed out by
Drinfeld \cite{D1}, it possesses a {\em pseudo-universal} R-matrix
$\RR(\lambda)\in(U_h(\hat\g')\otimes U_h(\hat\g'))((\lambda))$. The
$\lambda$ is related to the automorphism $\tau_\lambda$ defined in
\eqref{tau}. When using the homogeneous gradation $\RR(\lambda)$ is
a formal power series in $\lambda$.

When the pseudo-universal R-matrix is evaluated in the tensor
product of any two indecomposable finite-dimensional representations
$\rho_1$ and $\rho_2$ one obtains a numerical R-matrix
\begin{equation}\label{Rc}
 R^{12}(\lambda)=(\rho^1\otimes\rho^2)\RR(\lambda).
\end{equation}
The entries of these numerical R-matrices are rational functions of
the multiplicative spectral parameter $\lambda$ but when written in
terms of the additive spectral parameter $u = \log(\lambda)$ they
are trigonometric functions of $u$ and satisfy the Yang-Baxter
equation in the form given in \eqref{YB}. The matrix
$$
 \Rc^{12}(\lambda)= \sigma\circ R^{12}(\lambda)
$$
satisfies the intertwining relation
\begin{equation}\label{intertwining}
 \Rc^{12}(\lambda/\mu)\cdot(\rho^1_\lambda\otimes\rho^2_\mu)\Delta
 (x) = (\rho^2_\mu\otimes\rho^1_\lambda)(\Delta
 (x))\cdot \Rc^{12}(\lambda/\mu)
\end{equation}
for any $x\in U_h(\hat\g')$. It follows from the irreducibility of
the tensor product representations that these R-matrices satisfy the
Yang-Baxter equations
\begin{equation}\label{ryb}
 \begin{split}
 (\id\otimes\Rc^{23}(\mu/\nu))(\Rc^{13}(\lambda/\nu)\otimes\id)
 (\id\otimes\Rc^{12}(\lambda/\mu))\\
 =(\Rc^{12}(\lambda/\mu)\otimes\id)(\id\otimes\Rc^{13}(\lambda/\nu))
 (\Rc^{23}(\mu/\nu)\otimes\id)
 \end{split}
\end{equation}
or, graphically,
$$
\unitlength=1.4mm
\begin{picture}(70,26)
\thinlines \drawpath{14.0}{22.0}{34.0}{8.0}
\drawpath{34.0}{22.0}{14.0}{8.0} \drawpath{28.0}{22.0}{28.0}{8.0}
\drawpath{46.0}{22.0}{66.0}{8.0} \drawpath{66.0}{22.0}{46.0}{8.0}
\drawpath{52.0}{22.0}{52.0}{8.0} \drawthickdot{28.0}{18.0}
\drawthickdot{28.0}{12.0} \drawthickdot{52.0}{12.0}
\drawthickdot{52.0}{18.0} \drawthickdot{56.0}{15.0}
\drawthickdot{24.0}{15.0} \drawcenteredtext{40.0}{15.0}{$=$}
\drawcenteredtext{14.0}{25.0}{$V^3_\nu$}
\drawcenteredtext{21.0}{25.0}{$\otimes$}
\drawcenteredtext{28.0}{25.0}{$V^2_\mu$}
\drawcenteredtext{31.0}{25.0}{$\otimes$}
\drawcenteredtext{34.0}{25.0}{$V^1_\lambda$}
\drawcenteredtext{46.0}{25.0}{$V^3_\nu$}
\drawcenteredtext{49.0}{25.0}{$\otimes$}
\drawcenteredtext{52.0}{25.0}{$V^2_\mu$}
\drawcenteredtext{59.0}{25.0}{$\otimes$}
\drawcenteredtext{66.0}{25.0}{$V^1_\lambda$}
\drawcenteredtext{14.0}{5.0}{$V^1_\lambda$}
\drawcenteredtext{21.0}{5.0}{$\otimes$}
\drawcenteredtext{28.0}{5.0}{$V^2_\mu$}
\drawcenteredtext{31.0}{5.0}{$\otimes$}
\drawcenteredtext{34.0}{5.0}{$V^3_\nu$}
\drawcenteredtext{46.0}{5.0}{$V^1_\lambda$}
\drawcenteredtext{49.0}{5.0}{$\otimes$}
\drawcenteredtext{52.0}{5.0}{$V^2_\mu$}
\drawcenteredtext{59.0}{5.0}{$\otimes$}
\drawcenteredtext{66.0}{5.0}{$V^3_\nu$}
\end{picture}
$$

Explicit formulas for the pseudo-universal R-matrices were found
by Khoroshkin and Tolstoy. However these are difficult to evaluate
explicitly in specific representations so that in practice it is
easiest to find the numerical R-matrices $\Rc^{ab}(\lambda)$ by
solving the intertwining relation \eqref{intertwining}. It should
be stressed that solving the intertwining relation, which is a linear equation for the R-matrix,
is much easier than directly solving the Yang-Baxter equation, a cubic equation.

\section{Yangians}

As remarked by Drinfeld \cite{D2}, for untwisted $\hat\g$ the
quantum affine algebra $U_h(\hat\g')$ degenerates as $h\rightarrow
0$ into another quasi-pseudotriangular Hopf algebra, the {\em
Yangian} $Y(\g)$ \cite{D1}. It is associated with R-matrices which
are rational functions of the additive spectral parameter $u$. Its
representation ring coincides with that of $U_h(\hat\g)$.

Consider a general presentation of a Lie algebra $\g$, with
generators $I_a$ and structure constants $f_{abc}$, so that
$$[I_a,I_b]=f_{abc}I_c \,,\qquad \Delta(I_a) = I_a \otimes 1 + 1
\otimes I_a $$ (with summation over repeated indices). The Yangian
$Y(\g)$ is the algebra generated by these and a second set of
generators $J_a$ satisfying
$$ \left[ I_a , J_b \right] =   f_{abc} J_c \,, \qquad  \Delta(J_a)
= J_a \otimes 1 + 1 \otimes J_a + {1\over 2}f_{abc} I_c \otimes
I_b\,. $$ The requirement that $\Delta$ be a homomorphism imposes
further relations:
$$[J_a,[J_b,I_c]]-[I_a,[J_b,J_c]] =
\alpha_{abcdeg}\{I_d,I_e,I_g\} $$ and $$ [[J_a,J_b],[I_l,J_m]] +
[[J_l,J_m],[I_a,J_b]] = \left( \alpha_{abcdeg}f_{lmc} +
\alpha_{lmcdeg}f_{abc} \right) \left\{ I_d,I_e,J_g \right\} \,, $$
where
$$
 \alpha_{abcdeg}={1\over{24}}
 f_{adi} f_{bej}f_{cgk}f_{ijk}
 \hspace{0.1in},\hspace{0.15in}\{x_1,x_2,x_3\}=
 \sum_{i\neq j\neq k}x_{i}x_{j}x_{k} \,.
$$  When $\g={\mathfrak{sl}}_2$ the
first of these is trivial, while for $\g\neq {\mathfrak{sl}}_2$, the
first implies the second. The co-unit is $ \epsilon
(I_a)=\epsilon(J_a)=0$; the antipode is $s(I_a)=-I_a,\quad
s(J_a)=-J_a+{1\over 2}f_{abc}I_c I_b$. The Yangian may be obtained
from $U_h(\g^{(1)})$ by expanding in powers of $h$. For the precise
relationship see \cite{D1,MacK}. In the spin gradation, the
automorphism \eqref{tau} generated by $D$ descends to $Y(\g)$ as
$I_a \mapsto I_a,\; J_a\mapsto J_a + u I_a$.

There are two other realizations of $Y(\g)$. The first \cite{M}
defines $Y({\mathfrak{gl}}_n)$ directly from
$$R(u-v) T_1(u) T_2(v)=T_2(v) T_1(u) R(u-v),$$
where $T_1(u)=T(u)\otimes 1,\;T_2(v)=1\otimes T(v)$, and $$ T(u) =
\sum_{i,j=1}^n t_{ij}(u)\otimes e_{ij}\;,\qquad
t_{ij}(u)=\delta_{ij} + I_{ij}u^{-1} + J_{ij}u^{-2}+\ldots,$$
where $e_{ij}$ are the standard matrix units for
${\mathfrak{gl}}_n$. The rational R-matrix for the $n$-dimensional
representation of ${\mathfrak{gl}}_n$ is $$
 R(u-v) = 1 - {P\over u-v},\qquad {\rm where} \quad
P=\sum_{i,j=1}^n e_{ij}\otimes e_{ji}$$ is the transposition
operator. $Y({\mathfrak{gl}}_n)$ is then defined to be the algebra
generated by $I_{ij},J_{ij}$, and must be quotiented by the
`quantum determinant' at its centre to define
$Y({\mathfrak{sl}}_n)$. The coproduct takes a particularly simple
form,
$$ \Delta(t_{ij}(u)) = \sum_{k=1}^n t_{ik}(u)\otimes t_{kj}(u).$$

The third, Drinfeld's `new' realization of $Y(\g)$ \cite{D3}, we
do not give explicitly here, but we remark that it was in this
presentation that Drinfeld found a correspondence between certain
sets of polynomials and finite-dimensional irreducible
representations of $Y(\g)$, thus classifying these (although not
thereby deducing their dimension or constructing the action of
$Y(\g)$). As remarked earlier, the structure is as in sect.1.2:
$Y(\g)$ representations are in general $\g$-reducible, and there
is a set of $r$ fundamental $Y(\g)$-representations, containing
the fundamental $\g$-representations as components, from which all
other representations can be constructed.

\section{Origins in the Quantum Inverse Scattering Method}

Quantum affine algebras for general $\hat\g$ first appear in
\cite{D1,D2,J0,J}, but they have their origin in the `Quantum
Inverse Scattering Method' (QISM) of the St Petersburg school, and
the essential features of $U_h(\widehat{{\frak{sl}}_2})$ first
appear in \cite{KR}. In this section we explain how the quantization
of the Lax-pair description of affine Toda theory led to the
discovery of the $U_h(\hat\g)$ coproduct, commutation relations, and
R-matrix. We use the normalizations of \cite{J}, in which the $H_i$
are re-scaled so that the Cartan matrix $a_{ij}=\alpha_i.\alpha_j$
is symmetric.

We begin with the affine Toda field equations
$$
\partial^{\mu} \partial_{\mu} \phi_i = - {{m^2}\over{\beta}}
\sum_{j=1}^r \left( e^{\beta a_{ij} \phi_j} - n_i e^{\beta
\alpha_0. \alpha_j \phi_j} \right) \;,
$$
an integrable model in ${\mathbb R}^{1+1}$ of $r$ real scalar fields
$\phi_i(x,t)$ with a mass parameter $m$ and coupling constant
$\beta$. Equivalently, we may write $ \left[ \partial_x + L_x,
\partial_t + L_t \right] = 0 $ for the Lax pair
\begin{eqnarray*}
L_x(x,t) & = & {\beta \over 2} \sum_{i=1}^r H_i \partial_t\phi_i +
{m\over 2} \sum_{i,j=1}^r  e^{{\beta\over2}a_{ij}\phi_j} \left(
E^+_i + E^-_i \right) + {m\over 2}\sum_{j=1}^r
e^{{\beta\over2}a_{0j}
\phi_j} \left( \lambda E^+_0 + {1\over\lambda} E^-_0 \right)  \\
L_t(x,t) & = & {\beta \over 2} \sum_{i=1}^r H_i \partial_x \phi_i
 + {m\over 2} \sum_{i,j=1}^r
e^{{\beta\over2}a_{ij}\phi_j} \left( E^+_i - E^-_i \right) +
{m\over 2}\sum_{j=1}^r e^{{\beta\over2}a_{0j} \phi_j} \left(
\lambda E^+_0 - {1\over\lambda} E^-_0 \right)
\end{eqnarray*}
with arbitrary $\lambda\in{\mathbb C}$. The classical
integrability of the system is seen in the existence of
$r(\lambda,\lambda')$ such that
$$
\left\{ T(\lambda) \stackrel{\otimes}{,} T(\lambda')
 \right\} = \left[
r(\lambda,\lambda'), T(\lambda) \otimes T(\lambda') \right]\,,
$$ where $T(\lambda)=T(-\infty,\infty;\lambda)$
and $T(x,y;\lambda) = {\bf P} \exp \left( \int_x^y L(\xi;\lambda)
\,{d\xi} \right).$ Taking the trace of this relation gives an
infinity of charges in involution.

Quantization is problematic, owing to divergences in $T$. The QISM
regularizes these by putting the model on a lattice of spacing
$\Delta$, defining the lattice Lax operator to be
$$
 L_n(\lambda) =
 T((n-1/2)\Delta,(n+1/2)\Delta;\lambda) = {\bf P} \exp \left(
 \int_{(n-{1\over 2})\Delta}^{(n+{1\over 2})\Delta}
 L(\xi;\lambda) \,{d\xi} \right).
$$
The lattice monodromy matrix is then $T(\lambda)
=\lim_{l\rightarrow-\infty,m\rightarrow\infty}T_l^m$ where
$T_l^m=L_m L_{m-1}...L_{l+1}$, and its trace again yields an
infinity of commuting charges, provided that there exists a quantum
R-matrix $R(\lambda_1,\lambda_2)$ such that
\begin{equation}\label{rll}
 R(\lambda_1,\lambda_2) \, L_n^1(\lambda_1) \, L_n^2(\lambda_2)  =
 L_n^2(\lambda_2) \, L_n^1(\lambda_1) \, R(\lambda_1,\lambda_2),
\end{equation}
where $ L_n^1(\lambda_1)=L_n(\lambda_1) \otimes 1,
\;L_n^2(\lambda_2)=1 \otimes L_n(\lambda_2)$. That $R$ solves the
Yang-Baxter equation follows from the equivalence of the two ways of
intertwining $L_n(\lambda_1)\otimes L_n(\lambda_2)\otimes
L_n(\lambda_3)$ with $L_n(\lambda_3)\otimes L_n(\lambda_2)\otimes
L_n(\lambda_1)$.

To compute $L_n(\lambda)$, one uses the canonical, equal-time
commutation relations for the $\phi_i$ and $\dot{\phi}_i$. In
terms of the lattice fields
$$
 p_{i,n} = \int_{(n-{1\over 2})\Delta}^{(n+{1\over 2})\Delta}
  \dot{\phi}_i(x) \,{dx} \,,
 \hspace{0.15in}
  q_{i,n} = \int_{(n-{1\over 2})\Delta}^{(n+{1\over
 2})\Delta}\sum_j \,e^{{\beta\over2}a_{ij}\phi_j(x)} \,{dx} \,,
$$
the only non-trivial relation is $ \left[ p_{i,n} \, , \, q_{j,n}
\right] = {i\hbar\beta\over2}\delta_{ij} q_{j,n} \,, $ and one
finds
\begin{eqnarray*}
L_n(\lambda) & = & \exp \left( {\beta \over 2} \sum_i H_i p_{i,n}
\right) +
  \exp \left( {\beta \over 4} \sum_j H_j p_{j,n} \right) {m\over
2} \left[ \sum_{i} q_{i,n}\left( E^+_i
+ E^-_i \right) \right.\\[0.1in] \nonumber
& &\left. \hspace{0.4in} + \, \prod_i q^{-n_i}_{i,n} \left(
\lambda E^+_0 + {1\over\lambda} E^-_0 \right) \right]\exp
\left({\beta \over 4} \sum_j H_j p_{j,n} \right) + {
O}(\Delta^2) \,,
\end{eqnarray*}
the expression used by the St Petersburg school and by Jimbo. We now
make the replacement $E^\pm_i \mapsto q^{-H_i/4}E^\pm_i q^{H_i/4}$,
where $q=\exp(i\hbar\beta^2/2)$, and compute the ${ O}(\Delta)$
terms in \eqref{rll}, which reduce to
\begin{eqnarray*}
 R(z)(H_i \otimes 1 + 1\otimes H_i)  & = &
 ( H_i \otimes 1+1\otimes H_i)R(z) \\[0.1in]
 R(z) \left(
 E^{\pm}_i \otimes q^{-{H_i/2}}+ q^{{H_i/2}} \otimes E^{\pm}_i \right) & = &
 \left(
 q^{-{H_i/2}} \otimes E^\pm_i +E^{\pm}_i \otimes q^{{H_i/2}}\right) R(z)   \\[0.1in]
 R(z) \left(
 z^{\pm 1} E^{\pm}_0 \otimes q^{-{H_0/2}}+ q^{{H_0/2}} \otimes E^{\pm}_0 \right) & = & \left(
q^{-{H_0/2}} \otimes E^{\pm}_0+z^{\pm1} E^{\pm}_0 \otimes
 q^{{H_0/2}}  \right) R(z) \,,
\end{eqnarray*}
where $z=\lambda_1/\lambda_2$. We recognize in these the
$U_h(\hat\g)$ coproduct and thus the intertwining relations, in the
homogeneous gradation. These equations were solved for $R$ in
defining representations of non-exceptional $\g$ by Jimbo in
\cite{J}.

For $\hat\g=\widehat{{\mathfrak{sl}}_2}$, it was Kulish and
Reshetikhin \cite{KR} who first discovered that the requirement that
the coproduct must be an algebra homomorphism forces the replacement
of the commutation relations of $U(\widehat{\mathfrak{sl}}_2)$ by
those of $U_h(\widehat{\mathfrak{sl}}_2)$; more generally it
requires the replacement of $U(\hat\g)$ by $U_h(\hat\g)$.

\section{Affine quantum group symmetry and the exact S-matrix}

In the last section we saw the origins of $U_h(\hat\g)$ in the
`auxiliary' algebra introduced in the Lax pair. However, the quantum
affine algebras also play a second role, as a symmetry algebra. An
imaginary-coupled affine Toda field theory based on the affine
algebra $\hat\g^{\vee}$ possesses the quantum affine algebra
$U_h(\hat\g)$ as a symmetry algebra, where $\hat\g^{\vee}$ is the
Langland dual to $\hat\g$ (the algebra obtained by replacing roots
by coroots).

\begin{figure}[htb]
\unitlength=1.2mm
\begin{picture}(56,36)
\thinlines \drawpath{6.0}{28.0}{16.0}{14.0}
\drawpath{16.0}{28.0}{6.0}{14.0} \drawpath{36.0}{28.0}{42.0}{24.0}
\drawpath{42.0}{24.0}{48.0}{28.0} \drawpath{42.0}{24.0}{42.0}{18.0}
\drawpath{42.0}{18.0}{36.0}{14.0} \drawpath{42.0}{18.0}{48.0}{14.0}
\drawcenteredtext{6.0}{12.0}{$a$} \drawcenteredtext{16.0}{12.0}{$b$}
\drawlefttext{44.0}{20.0}{$c$} \drawcenteredtext{36.0}{12.0}{$a$}
\drawcenteredtext{48.0}{30.0}{$a$} \drawcenteredtext{6.0}{30.0}{$b$}
\drawcenteredtext{16.0}{30.0}{$a$}
\drawcenteredtext{48.0}{12.0}{$b$}
\drawcenteredtext{36.0}{30.0}{$b$}
\drawcenteredtext{42.0}{14.0}{$\theta^{ab}_c$}
\drawcenteredtext{11.0}{16.0}{$\theta$}
\drawpath{6.0}{28.0}{16.0}{14.0}\drawcenteredtext{10.0}{6.0}{a)}
\drawcenteredtext{42.0}{6.0}{b)}
\end{picture}
\caption{a) Graphical representation of a two-particle scattering
process described by the $S$-matrix $S_{ab}(\theta)$. b) At special
values $\theta_{ab}^c$ of the relative spectral parameter the two
particles of types $a$ and $b$ form a bound state of type
$c$.\label{sfig}}
\end{figure}
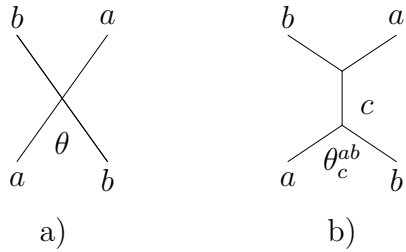

The solitonic particle states in affine Toda theories form
multiplets which transform in the fundamental representations of
the quantum affine algebra. Multi-particle states transform in
tensor product representations $V^a\otimes V^b$. The scattering of
two solitons of type $a$ and $b$ with relative rapidity $\theta$
is described by the $S$-matrix $S^{ab}(\theta) : V^a\otimes V^b
\rightarrow V^b\otimes V^a$, graphically represented in figure
\ref{sfig} a). It then follows from the symmetry that the
two-particle scattering matrix (S-matrix) for solitons must be
proportional to the intertwiner for these tensor product
representations, the $R$ matrix:
$$
 S^{ab}(\theta) = f^{ab}(\theta)\Rc^{ab}(\theta),
$$
with $\theta$ proportional to $u$, the additive spectral
parameter. The scalar prefactor $f^{ab}(\theta)$ is not determined
by the symmetry but is fixed by other requirements like unitarity,
crossing symmetry, and the bootstrap principle.

It turns out that the axiomatic properties of the R-matrices are in
perfect agreement with the axiomatic properties of the analytic
S-matrix. For example, crossing symmetry of the S-matrix,
graphically represented by
\begin{equation}
\unitlength=0.9mm \linethickness{0.4pt}
\begin{picture}(131.00,31.00)(5,25)
\put(12.00,50.00){\line(2,-3){14.00}}
\put(12.00,29.00){\line(2,3){14.00}}
\put(51.00,47.00){\line(3,-2){21.00}}
\put(51.00,33.00){\line(3,2){21.00}}
\put(19.00,35.00){\makebox(0,0)[ct]{$\theta$}}
\put(61.00,35.00){\makebox(0,0)[ct]{$i\pi-\theta$}}
\put(72.00,33.00){\line(0,-1){4.00}}
\put(72.00,27.00){\line(0,0){0.00}}
\put(51.00,47.00){\line(0,1){3.00}}
\put(51.00,41.00){\oval(16.00,16.00)[lb]}
\put(43.00,41.00){\line(0,1){9.00}}
\put(72.00,39.00){\oval(16.00,16.00)[rt]}
\put(80.00,39.00){\line(0,-1){10.00}}
\put(102.00,47.00){\line(3,-2){21.00}}
\put(102.00,33.00){\line(3,2){21.00}}
\put(112.00,35.00){\makebox(0,0)[ct]{$i\pi-\theta$}}
\put(102.00,39.00){\oval(16.00,16.00)[lt]}
\put(123.00,41.00){\oval(16.00,16.00)[rb]}
\put(131.00,41.00){\line(0,1){9.00}}
\put(94.00,39.00){\line(0,-1){10.00}}
\put(102.00,29.00){\line(0,1){4.00}}
\put(123.00,47.00){\line(0,1){3.00}}
\put(87.00,40.00){\makebox(0,0)[cc]{$=$}}
\put(34.00,40.00){\makebox(0,0)[cc]{$=$}}
\put(12.00,25.00){\makebox(0,0)[cb]{$a$}}
\put(26.00,25.00){\makebox(0,0)[cb]{$b$}}
\put(72.00,25.00){\makebox(0,0)[cb]{$a$}}
\put(80.00,25.00){\makebox(0,0)[cb]{$b$}}
\put(94.00,25.00){\makebox(0,0)[cb]{$a$}}
\put(102.00,25.00){\makebox(0,0)[cb]{$b$}}
\put(12.00,52.00){\makebox(0,0)[cb]{$b$}}
\put(26.00,52.00){\makebox(0,0)[cb]{$a$}}
\put(43.00,52.00){\makebox(0,0)[cb]{$b$}}
\put(51.00,52.00){\makebox(0,0)[cb]{$a$}}
\put(123.00,52.00){\makebox(0,0)[cb]{$b$}}
\put(131.00,52.00){\makebox(0,0)[cb]{$a$}}
\end{picture}
\end{equation}
is a consequence of the property of the universal R-matrix with
respect to the action of the antipode $S$,
$$
 (S\otimes 1)\RR = \RR^{-1}.
$$

An $S$-matrix will have poles at certain imaginary rapidities
$\theta^{ab}_c$ corresponding to the formation of virtual bound
states. This is graphically represented in figure \ref{sfig} b). The
location of the pole is determined by the masses of the three
particles involved,
$$
 m_c^2=m_a^2+m_b^2+2m_am_b\cos(i\theta^{ab}_c).
$$
At the bound state pole the $S$-matrix will project onto the
multiplet $V^c$. Thus the $\Rc$ matrix has to have this projection
property as well and indeed, this turns out to be the case. The
bootstrap principle, whereby the S-matrix for a bound state is
obtained from the S-matrices of the constituent particles,
\begin{equation}
\unitlength=0.9mm \linethickness{0.4pt}
\begin{picture}(96.00,42.00)(5,12)
\put(73.00,20.00){\line(2,3){10.00}}
\put(83.00,35.00){\line(2,-3){10.00}}
\put(83.00,35.00){\line(0,1){15.00}}
\put(20.00,20.00){\line(2,3){10.00}}
\put(30.00,35.00){\line(2,-3){10.00}}
\put(30.00,35.00){\line(0,1){15.00}}
\put(67.00,20.00){\line(5,3){28.00}}
\put(12.00,33.00){\line(5,3){28.00}}
\put(20.00,16.00){\makebox(0,0)[cb]{$a$}}
\put(40.00,16.00){\makebox(0,0)[cb]{$b$}}
\put(67.00,16.00){\makebox(0,0)[cb]{$d$}}
\put(73.00,16.00){\makebox(0,0)[cb]{$a$}}
\put(93.00,16.00){\makebox(0,0)[cb]{$b$}}
\put(30.00,52.00){\makebox(0,0)[cb]{$c$}}
\put(40.00,52.00){\makebox(0,0)[cb]{$d$}}
\put(83.00,52.00){\makebox(0,0)[cb]{$c$}}
\put(11.00,29.00){\makebox(0,0)[cb]{$d$}}
\put(96.00,39.00){\makebox(0,0)[cb]{$d$}}
\put(53.00,35.00){\makebox(0,0)[cc]{$=$}}
\end{picture}
\end{equation}
is a consequence of the property \eqref{DeltaR} of the universal R-matrix with
respect to the coproduct.

There is a famous no-go theorem (S. Coleman and J. Mandula,
Phys.Rev.\ 159 (1967) 1251) which states the ``impossibility of
combining space-time and internal symmetries in any but a trivial
way''. Affine quantum group symmetry circumvents this no-go theorem.
In fact, the derivation D is the infinitesimal two-dimensional
Lorentz boost generator and the other symmetry charges transform
non-trivially under these Lorentz transformations, see
\eqref{qrels}.

The non-cocommutative coproduct \eqref{coproduct} means that a
$U_h(\hat\g)$ symmetry generator, when acting on a two-soliton
state, acts differently on the left soliton than on the right
soliton. This is only possible because the generator is a non-local
symmetry charge -- that is, a charge which is obtained as the space
integral of the time component of a current which itself is a
non-local expression in terms of the fields of the theory.

Similarly, many nonlinear sigma models possess non-local charges
which form $Y(\g)$, and the construction proceeds similarly, now
utilising rational $R$-matrices, and with particle multiplets
forming fundamental representations of $Y(\g)$. In each case the
three-point couplings corresponding to the formation of bound
states, and thus the analogues for $U_h(\hat\g)$ and $Y(\g)$ of the
Clebsch-Gordan couplings, obey a rather beautiful geometric rule
originally deduced in simpler, purely elastic scattering models
\cite{CP2}.

More details about this topic can be found in \cite{Del2,MacK}.

\section{Integrable quantum spin chains}

Affine quantum groups provide an unlimited supply of integrable
quantum spin chains. From any R-matrix $R(\theta)$ for any tensor
product of finite-dimensional representations $W\otimes V$ one can
produce an integrable quantum system on the Hilbert space
$V^{\otimes n}$. This Hilbert space can then be interpreted as the
space of $n$ interacting spins. The space $W$ is an auxiliary space
required in the construction but not playing a role in the physics.

Given an arbitrary R-matrix $R(\theta)$ one defines the monodromy
matrix $T(\theta)\in\End(W\otimes V^{\otimes n})$ by
$$
 T(\theta) = R_{01}(\theta-\theta_1)R_{02}(\theta-\theta_2)
 \cdots R_{0n}(\theta-\theta_n)
$$
where, as usual, $R_{ij}$ is the R-matrix acting on the $i$-th and
$j$-th component of the tensor product space. The $\theta_i$ can be
chosen arbitrarily for convenience. Graphically the monodromy matrix
can be represented as
$$
\begin{picture}(180,60)
\unitlength=0.9mm \thinlines \drawpath{10.0}{16.0}{76.0}{16.0}
\drawpath{18.0}{22.0}{18.0}{10.0} \drawpath{28.0}{22.0}{28.0}{10.0}
\drawpath{38.0}{22.0}{38.0}{10.0} \drawpath{58.0}{22.0}{58.0}{10.0}
\drawpath{68.0}{22.0}{68.0}{10.0}
\drawcenteredtext{18.0}{6.0}{$V_1$} \drawpath{16.0}{6.0}{16.0}{6.0}
\drawcenteredtext{28.0}{6.0}{$V_2$}
\drawcenteredtext{38.0}{6.0}{$V_3$}
\drawcenteredtext{48.0}{6.0}{$\cdots$}
\drawcenteredtext{58.0}{6.0}{$V_{n-1}$}
\drawcenteredtext{68.0}{6.0}{$V_n$} \drawrighttext{8.0}{16.0}{$W$}
\end{picture}
$$
As a consequence of the Yang-Baxter equation satisfied by the
R-matrices the monodromy matrix satisfies
\begin{equation}\label{rtt}
 RTT=TTR.
\end{equation}
or, graphically,
$$
\unitlength=0.9mm
\begin{picture}(126,30)
\thinlines \drawpath{12.0}{14.0}{20.0}{22.0}
\drawpath{20.0}{22.0}{60.0}{22.0} \drawpath{12.0}{22.0}{20.0}{14.0}
\drawpath{20.0}{14.0}{60.0}{14.0} \drawpath{24.0}{26.0}{24.0}{10.0}
\drawpath{34.0}{26.0}{34.0}{10.0} \drawpath{56.0}{26.0}{56.0}{10.0}
\drawpath{110.0}{26.0}{110.0}{10.0}
\drawpath{114.0}{22.0}{122.0}{14.0}
\drawpath{122.0}{22.0}{114.0}{14.0}
\drawpath{70.0}{22.0}{114.0}{22.0} \drawpath{74.0}{26.0}{74.0}{10.0}
\drawpath{70.0}{14.0}{114.0}{14.0} \drawpath{84.0}{26.0}{84.0}{10.0}
\drawrighttext{10.0}{22.0}{$W$} \drawrighttext{10.0}{14.0}{$W'$}
\drawcenteredtext{24.0}{6.0}{$V_1$}
\drawcenteredtext{34.0}{6.0}{$V_2$}
\drawcenteredtext{56.0}{6.0}{$V_n$}
\drawcenteredtext{74.0}{6.0}{$V_1$}
\drawcenteredtext{84.0}{6.0}{$V_2$}
\drawcenteredtext{110.0}{6.0}{$V_n$}
\drawcenteredtext{44.0}{6.0}{$\cdots$}
\drawcenteredtext{96.0}{6.0}{$\cdots$}
\drawcenteredtext{64.0}{18.0}{$=$}
\end{picture}
$$
One defines the transfer matrix
$$
 \tau(\theta)=\tr_{W} T(\theta)
$$
which is now an operator on $V^{\otimes n}$, the Hilbert space of
the quantum spin chain. Due to \eqref{rtt} two transfer matrices
commute,
$$
 [\tau(\theta),\tau(\theta')]=0
$$
and thus the $\tau(\theta)$ can be seen as a generating function of
an infinite number of commuting charges, one of which will be chosen
as the Hamiltonian. This Hamiltonian can then be diagonalized using
the algebraic Bethe Ansatz.

One is usually interested in the thermodynamic limit where the
number of spins goes to infinity. In this limit, it has been
conjectured, the Hilbert space of the spin chain carries a certain
infinite-dimensional representation of the quantum affine algebra
and this has been used to solve the model algebraically, using
vertex operators \cite{JM}.

\section{Boundary quantum groups}

In applications to physical systems that have a boundary the
Yang-Baxter equation \eqref{YB} appears in conjunction with the
boundary Yang-Baxter equation, also known as the reflection
equation,
\begin{equation}\label{bYB}
 R_{12}(u-v)K_1(u)R_{21}(u+v)K_2(v)=K_2(v)R_{12}(u+v)K_1(u)R_{21}(u-v).
\end{equation}
The matrices $K$ are known as reflection matrices. This equation was
originally introduced by Cherednik to describe the reflection of
particles of a boundary in an integrable scattering theory and was
used by Sklyanin to construct integrable spin chains and quantum
field theories with boundaries.

Boundary quantum groups are certain co-ideal subalgebras of affine
quantum groups. They provide the algebraic structures underlying the
solutions of the boundary Yang-Baxter equation in the same way in
which affine quantum groups underlie the solutions of the ordinary
Yang-Baxter equation. Both allow one to find solutions of the
respective Yang-Baxter equation by solving a linear intertwining
relation. In the case without spectral parameters these algebras
appear in the theory of braided groups [See articles 28, 46].

For example the subalgebra $B_\epsilon(\hat\g)$ of $U_h(\hat\g')$
generated by
\begin{equation}
 Q_i=q_i^{H_i/2}(E^+_i+E^-_i)+\epsilon_i(q_i^{H_i}-1),\ \ i=0,\dots,r
\end{equation}
is a boundary quantum group for certain choices of the parameters
$\epsilon_i\in\C[[h]]$. It is a left coideal subalgebras of
$U_h(\hat\g')$ because
\begin{equation}
 \Delta(Q_i)=Q_i\otimes 1 + q_i^{H_i}\otimes Q_i\in
 U_h(\hat\g')\otimes B_\epsilon(\hat\g).
\end{equation}
Intertwiners $K(\lambda):V_{\eta\lambda}\rightarrow
V_{\eta/\lambda}$ for some constant $\eta$ satisfying
\begin{equation}
 K(\lambda) \rho_{\eta\lambda}(Q) = \rho_{\eta/\lambda}(Q)
 K(\lambda), \ \ \text{for all } Q\in B_\epsilon(\hat\g),
\end{equation}
provide solutions of the reflection equation in the form
\begin{equation}
 \begin{split}
 (\id\otimes K^2(\mu))\Rc^{12}(\lambda\mu)(\id\otimes
 K^1(\lambda))\Rc^{21}(\lambda/\mu)\\ =
 \Rc^{12}(\lambda/\mu)(\id\otimes K^1(\lambda))
 \Rc^{21}(\lambda\mu)(\id\otimes K^2(\mu)).
 \end{split}
\end{equation}
This can be expanded to the case where the boundary itself carries a
representation $W$ of $B_\epsilon(\hat\g)$. The boundary Yang-Baxter
equation can be represented graphically as
$$
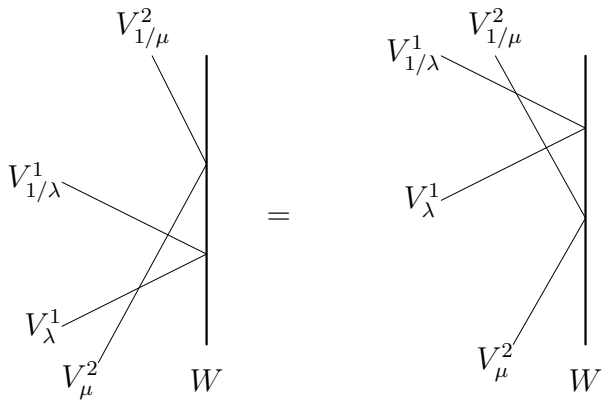
\unitlength=1.2mm
\begin{picture}(84,50)
\thicklines \drawpath{32.0}{42.0}{32.0}{10.0}
\drawpath{74.0}{42.0}{74.0}{10.0} \thinlines
\drawpath{32.0}{30.0}{26.0}{42.0} \drawpath{32.0}{30.0}{20.0}{8.0}
\drawpath{32.0}{20.0}{16.0}{28.0} \drawpath{32.0}{20.0}{16.0}{12.0}
\drawpath{74.0}{34.0}{58.0}{42.0} \drawpath{74.0}{34.0}{58.0}{26.0}
\drawpath{74.0}{24.0}{64.0}{42.0} \drawpath{74.0}{24.0}{66.0}{10.0}
\drawrighttext{16.0}{12.0}{$V^1_\lambda$}
\drawcenteredtext{18.0}{6.0}{$V^2_\mu$}
\drawrighttext{58.0}{26.0}{$V^1_\lambda$}
\drawcenteredtext{64.0}{8.0}{$V^2_\mu$}
\drawcenteredtext{32.0}{6.0}{$W$} \drawcenteredtext{74.0}{6.0}{$W$}
\drawrighttext{16.0}{28.0}{$V^1_{1/\lambda}$}
\drawcenteredtext{25.0}{45.0}{$V^2_{1/\mu}$}
\drawcenteredtext{64.0}{45.0}{$V^2_{1/\mu}$}
\drawrighttext{58.0}{42.0}{$V^1_{1/\lambda}$}
\drawcenteredtext{40.0}{24.0}{$=$}
\end{picture}
$$

Another example is provided by twisted Yangians where, when the
$I_a$ and $J_a$ are constructed as non-local charges in sigma
models, it is found that a boundary condition which preserves
integrability leaves only the subset
$$ I_i \qquad {\rm and} \qquad \tilde{J}_p = J_p + {1\over 4} f_{piq}(I_iI_q+I_qI_i)$$
conserved, where $i$ labels the $\h$-indices and $p,q$ the
${\mathfrak k}$-indices of a symmetric splitting $\g=\h+{\mathfrak
k}$. The algebra $Y(\g,\h)$ generated by the $I_i,\tilde{J}_p$ is,
like $B_\epsilon(\hat\g)$, a coideal subalgebra,
$\Delta(Y(\g,\h))\subset Y(\g)\otimes Y(\g,\h)$, and again yields
an intertwining relation for $K$-matrices. For
$\g={\mathfrak{sl}}_n$ and $\h={\mathfrak{so}}_n$ or
${\mathfrak{sp}}_{2n}$, $Y(\g,\h)$ is the {\em twisted Yangian}
described in \cite{M}.

All the constructions in earlier sections of this review have
analogs in the boundary setting. For more details see
\cite{Del,MacK}.


\baselineskip 14pt {\small
\begin{thebibliography}{10}

\bibitem{CP}
Chari V. and Pressley A.~N. (1994), Quantum Groups, CUP

\bibitem{CP2}
Chari V. and Pressley A.~N. (1996),Yangians, integrable quantum
systems and Dorey's rule, Commun. Math. Phys. 181, 265-302

\bibitem{Del}
Delius, G.~W. and MacKay, N.~J. (2003), Quantum Group Symmetry in
sine-Gordon and Affine Toda Field Theories on the Half-Line, Commun.
Math. Phys. 233, 173-190

\bibitem{Del2}
Delius, G.~W. (1995), Exact S-matrices with affine quantum group
symmetry, Nucl. Phys. B451, 445-465


\bibitem{D1}
Drinfeld V. (1985), Hopf algebras and the quantum Yang-Baxter
equation, Sov. Math. Dokl. 32, 254-258

\bibitem{D2}
 Drinfeld V. (1986), Quantum Groups, Proc. Int. Cong. Math.
(Berkeley) 798-820

\bibitem{D3}
Drinfeld V. (1988), A new realization of Yangians and quantized
affine algebras, Sov. Math. Dokl. 36, 212-216

\bibitem{J0}
Jimbo M. (1985), A $q$-difference analogue of $U(\g)$ and the
Yang-Baxter equation, Lett. Math. Phys. 10, 63-69

\bibitem{J}
Jimbo M. (1986), Quantum R-matrix for the generalized Toda system,
Commun. Math. Phys. 102, 537-547

\bibitem{JM}
Jimbo M. and Miwa T. (1995), Algebraic analysis of solvable
lattice models, AMS.

\bibitem{KR} Kulish P.~P. and Reshetikhin N.~Y. (1983),
Quantum linear problem for the sine-Gordon equation and higher
representations,
\newblock J. Sov. Math. 23, 2435.
\bibitem{MacK}
MacKay, N.~J. (2005), Introduction to Yangian symmetry in
integrable field theory, Int. J. Mod. Phys. to appear

\bibitem{M}
Molev A. (2003), Yangians and their applications, Handbook of
Algebra, Vol. 3, (M. Hazewinkel, Ed.), Elsevier, 907-959.





\end{thebibliography}

}
\end{document}